\documentclass[12pt]{article}
\usepackage{color,latexsym,amsfonts,amssymb,amsthm,amscd,amsmath,graphicx}
\newtheorem{lemma}{Lemma}
\newtheorem{corollary}{Corollary}

\newtheorem{theorem}{Theorem}
\newtheorem{proposition}{Proposition}

\def \min{{\rm min}}

\newcommand{\vy}{\mbox{\boldmath $y$}}

\newcommand{\vbeta}{\mbox{\boldmath $\beta$}}

\newcommand{\vepsilon}{\mbox{\boldmath $\epsilon$}}

\begin{document}
\baselineskip 8mm

\begin{center}
{\Large\bf  
Extended BIC for linear regression models with diverging number of relevant features and high or  ultra-high feature spaces \\}

\vspace{0.15in}

{\sc By}  SHAN LUO${}^{1}$ {\sc and}  ZEHUA CHEN${}^{*2}$\\
\vspace{0.15in}

${}^{1,2}$Department of Statistics and Applied Probability \\
National University of Singapore \\
3 Science Drive 2 \\
Singapore 117543 \\
Republic of Singapore \\ 
${}$ \\
Email: ${}^{1}$luoshan08@nus.edu.sg.  ${}^{2}$stachenz@nus.edu.sg, 
\vspace{1in}

Running title:   EBIC for linear models with diverging parameters
\end{center}

\newpage
\centerline{{\sc Summary}}

In many conventional scientific investigations with high or ultra-high dimensional feature spaces, the relevant features, though sparse, are large in number compared with classical statistical problems, and the magnitude of their effects tapers off.  It is reasonable to model the number of relevant features as a diverging sequence when sample size increases.  In this article, we investigate the properties  of the extended Bayes information criterion (EBIC) (Chen and Chen, 2008) for feature selection in linear regression models with diverging number of relevant features in  high or ultra-high dimensional feature spaces. The selection consistency of  the EBIC in this situation is established. The application of EBIC to feature selection is considered in a two-stage feature selection procedure.   Simulation studies are conducted to demonstrate the performance of the EBIC together with the two-stage feature selection procedure in finite sample cases.


\newpage

\section{ Introduction}

In the setting of a regression model, if the number of features (covariates) $p$ is of the polynomial order or exponential order of the sample size $n$, i.e., $p=O(n^\kappa)$  or  $p = O(\exp(n^\kappa))$, the feature space is referred to as a high-dimensional or ultra-high dimensional feature space.   Regression problems with high or ultra-high dimensional feature spaces arise in many important fields of scientific research such as genomics study, medical study, risk management, machine learning, etc..  Such problems are generally referred to as small-$n$-large-$p$ problems. In many small-$n$-large-$p$ problems the relevant (or causal, true, as referred by some other authors) features,  though sparse, are relatively large in number compared with classical statistical problems, and  their effects usually taper off to zero from the largest to the smallest.  To reflect the estimability of the feature effects,  it is reasonable to model the number of relevant features as a diverging sequence depending on the sample size.  \cite{Donoho00} and \cite{FanPeng04} are among the earliest papers dealing with diverging number of relevant features.  In this article, we consider model selection criteria for linear regression models with high or ultra-high feature space and diverging number of relevant features.

In general, there are two goals in model selection. The first one is to select models to do prediction and the focus is on prediction accuracy. The second one is to identify relevant features and the focus is on selection consistency.   In traditional model selection problems where the number of features under study is small,  these two goals might be addressed at the same time.  But, in small-$n$-large-$p$ problems,  the two goals need to be treated separately.  We concentrate on the second goal  in this article and refer to the problem as feature selection. 
A model selection criterion is crucial for feature selection.  The traditional model selection criteria such as Akaike's information criterion (AIC) \cite{Akaike73},  cross-validation (CV) \cite{Stone74},  generalized cross-validation (GCV) \cite{CravenWahba79}  and the Bayes information criterion (BIC) \cite{Schwarz78}   are not suitable for feature selection in small-$n$-large-$p$ problems.  The CV or GCV, which aims to minimize prediction errors, does not address the issue of selection consistency.  The AIC and BIC are overly liberal; that is, the criteria select far more features than the relevant ones, see,  \cite{BromanSpeed02, Siegmund04, BogdanDoergeGhosh04}. 
\cite{BogdanDoergeGhosh04}  proposed a modified BIC (mBIC) for the study  of genetic QTL mapping to address problems caused by too many features.  
\cite{ChenChen08} developed a family of extended Bayes information criteria (EBIC) for feature selection in small-$n$-large-$p$ problems. 
The family of EBIC is indexed by a parameter $\gamma$ in the range $[0, \ 1]$.  The original BIC is a special case of EBIC with $\gamma =0$. 
The mBIC is also a special case of EBIC in an asymptotic sense; that is, it is  asymptotically equivalent to the EBIC with $\gamma=1$. 
\cite{ChenChen08} considered the  case of high dimensional feature space with fixed number of relevant features.  They established the selection-consistency of EBIC when  $p=O(n^\kappa)$ and  $\gamma > 1- \frac{1}{2\kappa}$ for any $\kappa > 0$.  

 Model selection criterion for diverging number of relevant features in high or ultra-high dimensional feature space is still almost a void. \cite{WangLiLeng09} considered a BIC type criterion for diverging number of relevant features but their  criterion applies only when the dimension of the feature space is smaller than $n$, in fact,  they require $p/n^{\xi} <1$ for some $0<\xi<1$.  In this paper, we investigate the property of the EBIC when the number of relevant features diverges at the order $O(n^c)$ for some $0<c<1$ and  $p = O(n^\kappa)$ for any $\kappa$ or $p = O(\exp(n^\kappa))$ for some $0 < \kappa <1$.  We identify the conditions under which the EBIC remains selection consistent and provide the theoretical proof (Theorem 1).  
Since the seminal paper on LASSO \cite{Tibshirani96}, penalized likelihood methods with various penalty functions have been largely used for model selection, see, e.g., \cite{FanLi01,  MeinshausenYu09, Zhang10}.  It has been shown that if the penalty parameter in the penalized likelihood is properly chosen the penalized likelihood methods are selection consistent under certain conditions, see \cite{FanLi01, FanPeng04,  ZhaoYu06,  MeinshausenYu09, Wainwright06, KimChoiOh08}.  However, in practice, without a proper criterion for the selection of the penalty parameter (which corresponds to the selection of model), the selection consistency cannot be realized.  The commonly used criterion in the penalized likelihood methods, the CV, cannot be selection consistent in small-$n$-large-$p$ problems, as we have already pointed out in the previous paragraph.   In this paper, we  also consider the application of the EBIC for the selection of the penalty parameter in penalized likelihood methods.  Simulation studies are conducted to demonstrate the finite sample properties of the EBIC and the selection procedures. 

The remainder of the paper is arranged as follows. In \S 2, the  selection consistency of EBIC with diverging number of relevant features are established. In \S 3, a two-stage feature selection procedure with the application of the EBIC is described and discussed.   In \S 4,   simulation results are reported. Technical details and proofs are provided in the Appendix.

\section{Selection consistency of EBIC with diverging number of relevant features} 

We denote by $p_n$ the number of features under investigation to make its dependence on $n$ explicit. Let $(y_i, x_{i1}, \dots, x_{ip_n}),  i=1, \dots, n$,  be independent observations. We consider the following linear model
\begin{equation}
\label{lm}
y_{i} = \sum_{j=1}^{p_n} \beta_{nj} x_{ij} + \epsilon_{i}, \ i = 1, \dots, n,
\end{equation}
where $\epsilon_i$'s are i.i.d. with mean zero and variance $\sigma^2$. In matrix notation, (\ref{lm}) is expressed as
\[
 \vy_n = X_n \vbeta_n + \vepsilon_n, 
 \]  
where $\vbeta_n = (\beta_{n1}, \dots, \beta_{np_n})^T$,  $\vy_n = (y_1, \dots, y_n)^T$   and $X_n= (x_{ij})_{\stackrel{ i=1,\dots, n} {j=1,\dots,p_n}}$ .
 Here $p_n$ is either of a polynomial order or an exponential order of $n$, and $\vbeta_n$ is sparse, meaning that only a few of its components are non-zero. 
  
 We first introduce some notations.     Let $s_{0n} = \{ j: \beta_{nj} \neq 0, j \in \{ 1, \dots, p_n\} \}$.   
Let  $s$ be  any subset of $\{ 1, \dots, p_n\}$.  For convenience, we also refer to $s$ as a submodel. We denote by $X_n(s)$ the matrix composed of the columns of $X_n$ with  indices in $s$.  Similarly, $\vbeta_n(s)$ denotes the vector consisting of components of $\vbeta_n$ with indices in $s$.
Let $\nu(s)$ denote the number of components in $s$. Let $p_{0n} = \nu(s_{0n})$.  Let $H_n(s)$ be the projection matrix of $X_n(s)$, i.e., $H_n(s) = X_n(s) [ X_n(s)^T X_n(s)]^{-1} X_n(s)^T$.  Define
\[ \Delta_n(s) =  \| \mu_n - H_n(s) \mu_n\|_2^2, \]
where $\mu_n = E \vy_n = X_n(s_{0n}) \vbeta_n(s_{0n})$ and $\| \cdot \|_2$ is the $L_2$ norm. 

Let  ${\cal S}_j$ be the set of all combinations of $j$ indices in $\{ 1, \dots, p_n\}$.  Interchangeably we also call ${\cal S}_j$ the class of submodels consisting of $j$ features.
Let $\tau({\cal S}_j)$ be the size of ${\cal S}_j$; that is, $\tau({\cal S}_j)={p_n \choose j}$.

The family of EBIC proposed in \cite{ChenChen08} under model (\ref{lm}) is defined as
\[
\mbox{EBIC}_\gamma (s)
=
n   \ln \left( \frac{ \| \vy_n - H_n(s) \vy_n\|_2^2 }{ n }\right) 
   +  \nu (s) \ln n + 2 \gamma \ln \tau ({\cal S}_j), \ \  s \in {\cal S}_j,  \gamma \geq 0. 
\]
The family of EBIC is motivated from a Bayesian framework which gives rise to the BIC.  The BIC of a model $s$ is an approximation to the minus 2 log-transform of the posterior probability of $s$ while the prior probability on each model is equal.  With the equal prior probabilities,  the prior  probability on the submodel class $ {\cal S}_j$ is proportional to its size $  \tau ({\cal S}_j )$.  This makes BIC favor models with larger number of features in small-$n$-large-$p$ problems.  Instead of imposing an equal prior probability on each model, the EBIC imposes different prior probabilities on models in different submodel classes such that the prior probability on $ {\cal S}_j$ is proportional to  $  \tau ({\cal S}_j )^{-\gamma}$. The parameter $\gamma$ is determined such that the resultant EBIC is selection consistent.  In the case of high dimensional feature space, i.e., 
$p_n =O(n^{\kappa})$ for any $\kappa > 0$ , and a fixed number of relevant features, 
\cite{ChenChen08} showed that if $\gamma > 1 - 1/(2 \kappa)$ the EBIC is selection consistent. In the following, we deal with the case that the number of relevant features diverges and the feature space is high or ultra-high dimensional.    First we consider the following condition:

\noindent
{\bf Consistency Condition}:
{\it
\[
\lim_{n \to \infty}  \min 
\{  \frac{ \Delta_n(s)}{ p_{0n} \ln p_n} : s_0 \not\subset s,  \nu(s) \leq k_n \}  = \infty.
\]
where  $k_n = k p_{0n}$ for any fixed $k > 1$.  }

\noindent This condition is slightly different from what is called the asymptotic identifiability condition in \cite{ChenChen08}. 
The restriction $\nu(s) \leq k_n$ is imposed because in practice only the models with size comparable with and smaller than the true model will be considered.  Implicitly, the consistency condition requires that 
\begin{equation}
\label{CC1}
\sqrt{ \frac{ n}{ p_{0n}\ln p_n} }  \min\{ |\beta_{nj}|: j \in s_{0n}\}  \to \infty. 
\end{equation}
We now discuss a relationship between the consistency condition above and the well known sparse Reisz condition which is given as follows:
\begin{eqnarray*}
 0 < c_{\min}  & \leq & \min \{ \lambda_{\min} ( \frac{1}{n} X_n(s)^T X_n(s)) :  \nu(s) \leq k_n \} \\
& \leq & \max\{ \lambda_{\max} ( \frac{1}{n} X_n(s)^T X_n(s)) :  \nu(s) \leq k_n \}  \leq c_{\max} < \infty,
\end{eqnarray*}
where $\lambda_{\min}$ and $\lambda_{\max}$ denote the smallest and the largest eigenvalues respectively.
If $p_{0n}$ is fixed and hence so is $\{ \beta_{nj}: j \in s_{0n} \}$  then the sparse Reisz condition implies the consistency condition as shown in \cite{ChenChen08}.  If $p_{0n}$ diverges then the sparse Reisz condition together with (\ref{CC1}) imply the consistency condition.  When the number of relevant features diverges, conditions of the type (\ref{CC1}) are always imposed for selection consistency in penalized likelihood procedures, see \cite{ZhaoYu06,  Wainwright06, KimChoiOh08}. As the following proposition implies, the sparse Reisz condition together with (\ref{CC1}) are stronger assumptions than the consistency condition. 
\begin{proposition}
 Assume $s_{0n} = \{ 1, 2, \dots, p_{0n}\}$.
Let $s_{-k}$ be the set with the $k$th element of $s_{0n}$ removed.
Let $ k(s) = s_{-k} \cup s$. If  (\ref{CC1}) is satisfied and
\begin{equation}
\label{SC}
\lim_{n \to \infty}
\min_{s:  \nu(s) \leq k_n, s_0 \not\subset s }
 \frac{   \max_{k} \{   \| [ I - H_n(k(s))] X_n(\{k\}) \|  \}}{ p_{0n} \ln p_n } = \infty
\end{equation}
then the consistency  condition holds.
\end{proposition}
\noindent The above proposition is similar to a result in \cite{ChenChen08} which deals with a high dimensional feature space and a fixed number of relevant features.   The same  as in \cite{ChenChen08},   examples can be constructed such that (\ref{SC}) holds but the sparse Reisz condition does not hold.  

Condition (\ref{CC1}) determines the divergence pattern of $(n, p_{0n}, p_n)$ and the constraint on $\beta_{nj}$. Now consider the high and ultra-high dimensional feature spaces separately.   If $p_n = O(n^{\kappa})$ for any fixed $\kappa > 0$ and  $p_{0n} = n^{c}$ for some $0<c <\kappa$,  (\ref{CC1}) reduces to  
$ \frac{ n^{1-c}}{\ln n}   \min\{ |\beta_{nj}^2|: j \in s_{0n}\}  \to \infty. $ The induced constraint on $\beta_{nj}$ is that $\min\{ |\beta_{nj}^2|: j \in s_{0n}\}$ must have a magnitude larger than $O(  n^{-(1-c)} )$. Let $b$ be any number bigger than $c$. Then the following  provides a consistency pattern:  $(n, p_{0n}, p_n) = ( n, O(n^c), O(n^\kappa))$, $\min\{ |\beta_{nj}|: j \in s_{0n}\} = O(n^{-(1-b)/2})$,  $0<c<\kappa$, $c<b<1$.
If $p_n = O(\exp(n^\kappa))$ and  $p_{0n} = n^{c}$ then, by the same argument,  (\ref{CC1}) induces the following consistency pattern: $(n, p_{0n}, p_n) = ( n, O(n^c), O(\exp(n^\kappa) ) )$, $\min\{ |\beta_{nj}|: j \in s_{0n}\} = O(n^{-(1-b)/2})$,  $0<c , \kappa <1 $, $c+\kappa < b<1$.

We now state the main result on the selection consistency of the EBIC with diverging number of relevant features in high or ultra-high dimensional feature spaces.

\begin{theorem} 
Assume  model (\ref{lm}) and the consistency condition.  In  addition, assume that $p_{0n} \ln p_n = o(n)$,  $\ln p_{0n}/\ln p_n \to \delta \geq 0$. 
Let $k_n = kp_{0n}$ for any constant $k>1$.  Then 
\[
P \{ \min_{s: \nu(s) \leq k_n }   \mbox{EBIC}_{\gamma} (s) 
> \mbox{EBIC}_{\gamma}(s_{0n}) \} \to 1,
\]
 if $\gamma > \frac{1+\delta}{1-\delta} - \frac{\ln n }{2 (1-\delta) \ln p_n}$. 
\end{theorem}

\noindent The following are immediate corollaries of Theorem 1.

\begin{corollary}
 If  $p_n = O(n^\kappa)$ for any constant $\kappa >0 $,  $p_{0n} = p_0$  is fixed,  the EBIC is selection consistent with  
 $\gamma > 1 - \frac{\ln n }{2 \ln p_n} = 1-\frac{1}{2\kappa} $ among all models $s$ with  $\nu(s) \leq k_n $.
\end{corollary}

\begin{corollary}
 If  $p_n = O(n^\kappa)$ for any constant $\kappa >0 $,  $p_{0n} = O(n^{c})$,  $\min\{ |\beta_{nj}|: j \in s_{0n}\} = O(n^{-(1-b)/2})$,  $0<c<\kappa$, $c<b<1$,    then the EBIC is selection consistent with  
 $\gamma > \frac{\kappa + c -0.5}{\kappa - c} $ among all models $s$ with  $\nu(s) \leq k_n $.
\end{corollary}

\begin{corollary}
 If  $p_n = O(\exp(n^\kappa))$ for  $0< \kappa <1$, $p_{0n} = O(n^{c})$,  $\min\{ |\beta_{nj}|: j \in s_{0n}\} = O(n^{-(1-b)/2})$, $0<c , \kappa <1 $, $c+\kappa < b<1$, 
the EBIC is selection consistent with  
 $\gamma > 1 - \frac{\ln n }{2 \ln p_n}$ among all models $s$ with  $\nu(s) \leq k_n $.
\end{corollary}

\noindent The following  lemmas are needed in the proof of Theorem 1. 

\begin{lemma}
\label{lemma1} If  $\dfrac{\ln j}{\ln p}\rightarrow \delta \;\mbox{as}\;p\rightarrow +\infty,$ we have
\[ \ln (\dfrac{p!}{j!(p-j)!})=j\ln  p (1-\delta)(1+o(1)).
\]
\end{lemma}

\begin{lemma}
\label{lemma2}  Let $\chi^2_k$ denote a $\chi^2$ random variable with degrees of freedom $k$. If  $m \rightarrow
+\infty$ and $\dfrac{K}{m}\rightarrow 0$ then 
\[
P(\chi^2_{k}\geq m)=
\dfrac{1}{\Gamma(k/2)}(m/2)^{k/2-1}e^{-m/2}(1+o(1)),
\]
uniformly for all $k \leq K$.
\end{lemma}

\noindent The proofs of   Lemmas 1 and 2 and Theorem 1 are given in the Appendix.

\section{Application of EBIC in feature selection procedures}

In this section. we consider the application of  EBIC for  choosing tuning parameters in feature selection procedures using penalized likelihood methods.  When the dimension of the feature space is high or ultra-high,  a natural first step in feature selection is  to reduce the dimensionality of the feature space by some screening procedure and then to apply the penalized likelihood method with the reduced feature space.  This has become a well-accepted strategy for feature selection with high or ultra-high feature space, see, e.g.,  \cite{FanLv08, WassermanBoeder09, ChenChen09}. In the following, we describe a general feature selection procedure of this nature where EBIC is used to choose the penalty parameter in the penalized likelihood.   
\begin{description}
\item [Screening stage:]  Let ${\cal F}_n$  denote the set of all the features.  This stage  screens out obviously irrelevant features by a screening procedure and reduces  ${\cal F}_n$  to a set ${\cal S}_n^*$ with dimension smaller than $n$.   The screening procedure we recommend is as follows. First using the sure independence screening (SIS) advocated in  \cite{FanLv08} to reduce the dimension of ${\cal F}_n$ to a low power order of $n$, say $n^{3/2}$, then using LASSO by choosing an appropriate penalty  parameter to further reduce the dimension below $n$.   
\item [Selection stage:]  Select features by optimizing  a penalized log likelihood of the form
\[ l_{n, \lambda} (X({\cal S}_n^*), \vbeta({\cal S}_n^*) ) = -2 \ln L(X({\cal S}_n^*), \vbeta({\cal S}_n^*) )  + \sum_{j \in {\cal S}_n^*} p_{\lambda}(|\beta_j|), \]
where $L(X({\cal S}_n^*), \vbeta({\cal S}_n^*) )$ is the likelihood function of the model with all features in ${\cal S}_n^*$ , $p_{\lambda}(\cdot)$ is a penalty function and $\lambda$ is the penalty parameter. An appropriate penalty function to use is the SCAD penalty proposed in \cite{FanLi01}.  The $\lambda$ is chosen by EBIC as follows.  For each $\lambda$, let $s_{n \lambda}$ be the set of features with non-zero coefficient when $l_{n, \lambda} (X({\cal S}_n^*), \vbeta({\cal S}_n^*) ) $ is minimized.  Compute 
\[ \mbox{EBIC}_{\gamma} (\lambda) = -2 \ln L(X(s_{n \lambda}), \hat{\vbeta}(s_{n \lambda}) )  + \nu(s_{n \lambda}) \ln n  + 2 \gamma \ln {p_n \choose \nu(s_{n \lambda}) },\]
where $\hat{\vbeta}(s_{n \lambda}) $ is the maximum likelihood estimate (without penalty) of $\vbeta(s_{n \lambda}) $ and $\gamma$ is taken to be $1- \frac{\ln n}{ C \ln p_n}$ for some $C >2$. 
Let $\lambda^*$ be the one which attains the minimum $\mbox{EBIC}_{\gamma}$. The set $s_{n \lambda^*}$ is taken as the set of selected features.
\end{description}
We shortly discuss the properties of the above feature selection procedure in the following. 
For a screening procedure, if  $P({\cal S}_n^* \subset {\cal F}_n) \to 1$, as $n$ goes to infinity,  the screening procedure is said to have the property of sure screening, see \cite{FanLv08}.   For a penalized likelihood function of the above type, if  there is $\lambda_n$ such that
$ P (s_{n \lambda_n}   = s_{0n}) \to 1, $ the penalized likelihood is said to have an oracle property (in a narrower sense).   
If the screening procedure in the screening stage has the property of sure screening,  the reduced feature space ${\cal S}^*_n$  will contain all the relevant features with probability converging to 1 as $n$ goes to infinity.  If the penalized likelihood has the oracle property with the reduced feature space, there will be a $\lambda$ value such that its corresponding set $s_{n \lambda}$ is the same as $s_{0n}$, the true set of relevant features,   in the selection stage when ${\cal S}^*_n$ contains all the relevant features.  Then the selection consistency of EBIC will guarantee that the true set of relevant feature is selected.  Thus the  feature selection procedure will be selection consistent if the conditions required by the sure screening property of the screening procedure, the oracle property of the penalized likelihood and the selection consistency of EBIC are met simultaneously.

Fan and Lv \cite{FanLv08} showed that, under certain conditions (conditions 1-4 in section 5 of their paper),  the SIS has the sure screening property if the dimension of the feature space is reduced to an order $O(n^{1-\theta})$ for some $\theta > 0$.  If the tuning parameter in LASSO is chosen such that the number of non-zero coefficients  is large enough (smaller than $n$),  the LASSO procedure can retain all the true features almost surely as $n$ goes to infinity, see \cite{ChenChen09}.  Kim {\it et al.} 
  \cite{KimChoiOh08} considered the  SCAD with diverging number of relevant features under the following conditions.
C1:  There are $0 \leq c < b \leq 1$ and $M_1>0$ such that 
	$  p_{0n} = O(n^{c})  \   \mbox{and} \  n^{(1-b)/2} \min \{|\beta_{nj}|:  j \in s_{0n} \} \geq M_1. $
C2: There exists $M_2>0$ such that
$  n^{-1} (X_n(\{j \})^T X_n(\{j \}) \leq M_2,  \  \mbox{for any}  \  j .$ 
C3: There exists $M_3 > 0$ such that 
$ \lambda_{\min} ( n^{-1} X_n(s_{n0})^T X_n (s_{n0}) ) \geq M_3 , $
where $ \lambda_{\min}$ denotes the smallest eigenvalue. 
C4:  $p_n \leq n$ and the eigenvalues of  $n^{-1} X_n^T X_n $ are uniformly bounded from both below and above. 
They showed that under the above conditions the oracle property of the SCAD holds. The condition  $ n^{(1-b)/2} \min \{|\beta_j|:  j \in s_{0n} \} \geq M_1 $ implies that  $n^{(1-c-\kappa)/2} \min \{|\beta_j|:  j \in s_{0n} \}  \to \infty$ for some $\kappa < b-c$. 
If  C4 is replaced by C$4^{'}$: $\tau({\cal S}_n^*) < n $ and  the eigenvalues of  $n^{-1} X_n(s_{n})^T X_n (s_{n})$  for any $s_n \subset {\cal S}_n^* $ are uniformly bounded from both below and above,  then together with C1-C3 the oracle property of the SCAD penalized likelihood in the selection stage will be guaranteed. 
 Therefore,     suppose that conditions 1-4 in \cite{FanLv08}, C1-C3, C$4^{'}$ and the consistency condition hold,  then the two-stage procedure described above  is selection consistent.    The reason we recommend a two-step screening procedure is that if only SIS is used to reduce the dimensionality below $n$  condition C$4^{'}$ might not hold because SIS does not reduce the level of the spurious correlations in the original feature space.  On the other hand, LASSO does reduce the level of the spurious correlations since it tends to select only one of the highly correlated features, see \cite{ZouHastie05}, but due to the capacity of the computing facilities it might not be able to handle ultra-high dimensional feature space. When the two steps are combined it is more likely that C$4^{'}$ will be satisfied while the sure screening property is retained.  In fact, the conditions in  \cite{FanLv08} for the sure screening property can be much relaxed when the dimensionality is only reduced to a power order of $n$ higher than $n$.  The performance of the feature selection procedure described in this section is investigated in simulation studies which are presented in the next section.

\section{Simulation studies }

The purpose of the simulation studies is to investigate the applicability of EBIC in feature selection procedures and to investigate whether or not the asymptotic property of selection consistency can be realized in finite sample situations. To this end, the two-stage feature selection procedure discussed in \S 3 is considered in the simulation studies.   The R package \verb+plus+   \cite{Zhang07} is used for the computation.  We are mainly concerned about the selection consistency of  the EBIC in the consistent range of $\gamma$.  We take $\gamma$ slightly bigger than $ 1- \frac{\ln n} { 2 \ln p_n}$ (in the simulation we take  $\gamma =1- \frac{\ln n} { 4 \ln p_n}$)  for  demonstrating the performance of the EBIC in  finite sample situations.  We also consider $\gamma = 0$, which corresponds to the original BIC,  and $\gamma = 1$, which corresponds to an asymptotic form of the mBIC proposed in \cite{BogdanDoergeGhosh04}.  Throughout the simulation studies, 
$\tau({\cal S}_n^*)$ is taken to be $0.5n$. 

We take the  divergence pattern as $(n,   p_{0n}, p_n) = (n,  c[n^{0.325}], [\exp(n^{0.35})])$ for $n=100, 200, 500$ and $1,000$, and $c=1$ and 2, which results in the table below:
\begin{center}
\begin{tabular}{c|cccc} \hline
 $n$ &  100  &200 & 500  &1,000 \\
   $p_n$ &   150&  595 &6,655 &74,622 \\ 
$p_{0n}(c=1)$ &    4  &  6   & 8   &  9 \\
$p_{0n}(c=2)$  &    8  &  12   & 16   &  18 \\\hline
\end{tabular}
\end{center}
For $ j \in s_{0n}$ the parameter $\beta_{nj}$ is independently generated as 
$\beta_{nj} = (-1)^u$ $(n^{-0.1625}+|z|)$ where
$ u\sim Bernoulli(0.4)$ and $z$ is a normal random variable
with mean $0$ and satisfies $P(|z|\geq 0.1)=0.25$.  This ensures, roughly, $\min\{|\beta_{nj}|: j \in s_{0n}\} = O(n^{-0.1625})$.
The error variance $\sigma^2$ is determined by setting the following ratio to certain values when $n=100$ and kept unchanged for other $n$'s:
\[ h =\dfrac{ E ({\vbeta^*}^T\Sigma \vbeta^* ) }{E({\vbeta^*}^T\Sigma \vbeta^*)+\sigma^2},\]
where $\Sigma$ is  the covariance matrix of the predictors and  the expectation is with respect to the generating distribution of $\vbeta^*$.   This ratio mimics what is  called the heritability in broad sense in genetic studies.   We considered  $h = 0.4, 0.6$ and $0.8$.  For each simulation setting, 
 $200$ data sets are generated and analyzed.  
The following three correlation structures are considered for the covariates:

{\it Structure I:  Power decay correlation.} The covariates are generated as a series of  normally distributed random variables with 
mean $0$ and correlation coefficient $\rho_{ij} = 0.5^{|i-j|}$.

{\it Structure II:  Diagonal block design with equal pairwise correlation. } 
The  covariance matrix  is a diagonal block matrix. Each block except the last one is of dimension $50 \times 50$. The variances in the blocks are all equal to 1 and the off-diagonal correlations are all equal to $\rho=0.5$.

{\it Structure III:  Diagonal block design with  uniformly distributed eigenvalues. } 
Unlike the diagonal block matrix in Structure 2, each block is first generated such that its smallest eigenvalue is 1,  largest eigenvalue is 50 and  other eigenvalues are uniformly distributed over $(1, 50)$,  and then it is converted into a correlation matrix. 

The finite sample performance of the EBIC is assessed by the positive discovery rate (PDR) and false discovery rate (FDR) defined as follows:
\[ \mbox{PDR}_n =\frac{\nu(s_{n\lambda^*} \cap s_{0n})}{\nu(s_{0n})}, \ \ \mbox{FDR}_n =\frac{\nu(s_{n\lambda^*}\backslash s_{0n})}{\nu(s_{n\lambda^*})},\]
where $s_{n\lambda^*}$ is the set of features selected in the selection stage of the two-stage procedure.  The asymptotic property of selection consistency is equivalent to 
\[ \lim_{n\to \infty }  \mbox{PDR}_n  = 1  \ \ \mbox{and} \ \ \lim_{n\to \infty }  \mbox{FDR}_n  = 0, \]
in probability.   

The simulated $\mbox{PDR}_n$  and $\mbox{FDR}_n$  averaged over  200 replicates for each setting are  reported in Table~\ref{table1}, \ref{table2} and \ref{table3} respectively for correlation stucture I, II and III.  In the tables,   $\gamma_{\mbox{\tiny BIC}} =0$   corresponding to BIC,  $\gamma_{\mbox{\tiny SC}} = 1- \ln n/(4 \ln p_n)$ corresponding to a value in the selection consistent range of $\gamma$ and 
$\gamma_{\mbox{{\tiny  mBIC}}} =1 $ corresponding to mBIC.

The following points are manifest in Tables 1, 2 and 3.  (i) The finite sample performance of the EBIC closely matches its asymptotic property. That is, under all the three correlation structures, for the procedure with  $\mbox{EBIC}_{\gamma_{\mbox{\tiny SC}}}$, the $\mbox{PDR}_n$   and the $\mbox{FDR}_n$ approach rapidly to 1 and 0 respectively,  as $n$ increases from 100 to 1000, at all the three $h$ levels.   (ii) The  BIC does  not appear to be selection consistent.  Under all the settings,  the $\mbox{FDR}_n$ of the procedure with BIC does not reduce as $n$ increases, it is in fact quite the opposite.   (iii) In general,  the $\mbox{PDR}_n$ of the procedure with BIC is  higher because  it always selects much more features.  But, as $n$ gets large,  the $\mbox{PDR}_n$ of   $\mbox{EBIC}_{\gamma_{\mbox{\tiny SC}}}$ quickly becomes comparable with that of the BIC.  (iv) For large $n$, the mBIC is comparable with $\mbox{EBIC}_{\gamma_{\mbox{\tiny SC}}}$, which  reflects the fact that it is also  selection consistent since $\gamma_{\mbox{\tiny mBIC}} = 1$ is in the consistency range of EBIC.  But for small $n$,  it loses certain power while overly controlling $\mbox{FDR}_n$.

\appendix

\section{Appendix}\label{app}

\subsection{Proof of Lemma 1:}

\begin{proof}
\noindent Write 
\[ \dfrac{p!}{j!(p-j)!} = \frac{ p (p-1) \cdots (p-j+1)  }{j!} =  \dfrac{ p^j   \left(1-\frac{1}{p}\right) \cdots  \left(1- \frac{j-1}{p}\right)} { j!}.\]
  Note that 
\[
 \left(1-\frac{j-1}{p}\right)^{j-1} <  \left(1-\frac{1}{p}\right) \cdots  \left(1- \frac{j-1}{p}\right) <  \left( 1 - \frac{1}{p}\right)^{j-1},  
\]
and, see \cite{Robbins55},  that 
\[
\sqrt{2\pi}j^{j+1/2}e^{-j+1/(12j+1)}<j!<\sqrt{2\pi}j^{j+1/2}e^{-j+1/(12j)}.  
\]
We now have 
\begin{equation} 
\label{eq1}
\begin{split}
\ln (\dfrac{p!}{j!(p-j)!})\leq &  j\ln p + (j-1) \ln (1-1/p)   - (j+1/2) \ln j +j - \frac{1}{12j+1} - \ln \sqrt{2 \pi} \\
               \leq &  j\ln p  - (j+1/2) \ln j + j  
                =  j\ln p [1  -  \frac{(j+1/2) \ln j}{j \ln p} + \frac{1}{\ln p} ] \\
               = &  j\ln p (1-\delta) (1+o(1)).
\end{split}
\end{equation}
and 
\begin{equation}
\label{eq2}
\begin{split}
\ln (\dfrac{p!}{j!(p-j)!})\geq&  j\ln p + (j-1) \ln (1-\frac{j-1}{p})   \\
                   & - (j+1/2) \ln j +j - \frac{1}{12j} - \ln \sqrt{2 \pi} \\
           \geq &  j\ln p + (j-1) \ln \left(1-\frac{j-1}{p}\right)  
                    - (j+1/2) \ln j - \ln \sqrt{2 \pi} \\
                = & j\ln p\left(1  + \frac{(j-1) \ln \left(1-\frac{j-1}{p}\right)  }{j \ln p}
                    - \frac{(j+1/2)  \ln j}{j \ln p}  - \frac{\ln \sqrt{2 \pi} }{j \ln p}\right) \\
                =& j\ln p (1-\delta)(1+o(1)).  
\end{split}
\end{equation}
Lemma1 follows from (\ref{eq1}) and (\ref{eq2}).
\end{proof}

\subsection{Proof of Lemma 2}

\begin{proof}

\noindent Denote $\bar{F}_{k}(m)=P(\chi^2_{k}\geq m)$. By integration by parts, we obtain 
\[
\bar{F}_{k}(m)= \dfrac{1}{2^{k/2}\Gamma(k/2)}\int_{m}^{+\infty}x^{k/2-1}e^{-x/2} dx
=\dfrac{1}{\Gamma(k/2)}(m/2)^{k/2-1}e^{-m/2}+\tilde{F}_{k-2}(m).
\]
If $k$ is even, 
\[ \bar{F}_{k}(m)
= \dfrac{1}{\Gamma(k/2)}(m/2)^{k/2-1}e^{-m/2}[1+\sum\limits_{i=1}^{(k-2)/2}(\dfrac{(k/2-1)\cdots(k/2-i)}{(m/2)^i})].
\]
If $k$ is odd, 
\[ \bar{F}_k(m) = 
\dfrac{1}{\Gamma(k/2)}(m/2)^{k/2-1}e^{-m/2}[1+\sum\limits_{i=1}^{(k-3)/2}(\dfrac{(k/2-1)\cdots(k/2-i)}{(m/2)^i})]+\bar{F}_1(m),
\]
where 
\[ \bar{F}_1(m)=P(\chi^2_1\geq m)\approx
2\dfrac{\exp(-m/2)}{\sqrt{2\pi m}} =  \dfrac{1}{\Gamma(k/2)}  (m/2)^{k/2-1}e^{-m/2} \dfrac{2\Gamma(k/2)}{\sqrt{2\pi}(m/2)^{(k-1)/2}}\] when $m \rightarrow
+\infty. $  We can write 
\[ \bar{F}_k(m) = 
\dfrac{1}{\Gamma(k/2)}(m/2)^{k/2-1}e^{-m/2}[1+ R(k,m)]. 
\]
It is straightforward to see that $R(k,m) \leq R(K,m) \to 0$ when $m\to +\infty.$ 
\end{proof}

\subsection{Proof of Theorem 1}

\begin{proof}

\noindent Let $s$ be any submodel.  Decompose $\mbox{{\rm EBIC}}_{\gamma}(s)-\mbox{{\rm EBIC}}_{\gamma}(s_{0n})$ 
as follows:

\begin{equation}
\begin{split}
  & \mbox{{\rm EBIC}}_{\gamma}(s) - \mbox{{\rm EBIC}}_{\gamma}(s_{0n}) \\ \label{peq1}
=&n\ln \dfrac{y_n^T\{I_n-H_n(s)\}y_n}{y_n^T\{I_n-H_n(s_{0n})\}y_n}+(\nu(s) -p_{0n})\ln
n+2\gamma(\ln \tau(\mathcal{S}_j)-\ln \tau(\mathcal{S}_{p_{0n}}))\\
=& T_1 + T_2,  \ \mbox{say},  
\end{split}
\end{equation}
where 
\begin{equation}
\label{peq2}
\begin{split}
 T_1 = & n \ln \dfrac{y_n^T[I_n-H_n(s)]y_n}{y_n^T[I_n-H_n(s_{0n})]y_n}
                =  n \ln \dfrac{y_n^T[I_n-H_n(s)]y_n}{\epsilon_n^T[I_n-H_n(s_{0n})]\epsilon_n}\\
        = & n\ln \left\{1+ \dfrac{y_n^T[I_n-H_n(s)]y_n-\epsilon_n^T[I_n-H_n(s_{0n})]\epsilon_n}
              {\epsilon_n^T[I_n-H_n(s_{0n})]\epsilon_n}\right\} \\
 T_2 = & (\nu(s)-p_{0n})\ln n+2\gamma(\ln \tau(\mathcal{S}_{\nu(s)} )-\ln \tau(\mathcal{S}_{p_{0n}})).
\end{split}
\end{equation}

\vspace{0.15in}
\noindent {\it Case I:} $s_{0n} \not\subset s$.  
\vspace{0.15in}

Without loss of generality, assume $\sigma^2 =1$. We can write 
\begin{equation}
\label{peq3}
\epsilon_n^T\{I_n-H_n(s_{0n})\}\epsilon_n=\sum\limits_{i=1}^{n-p_{0n}}Z_i^2=(n-p_{0n})(1+o_p(1))=n(1+o_p(1)),
\end{equation}
where $Z_i$'s are i.i.d. standard normal variables, since $H_n(s_{0n})$ is a projection matrix with rank $p_{0n}$.  
We have  
\[
\begin{split}
&y_n^T[I_n-H_n(s)]y_n-\epsilon_n^T[ I_n-H_n(s_{0n}) ]\epsilon_n \\
=&\Delta_n(s)+2\mu_n^T[ I_n-H_n(s)]\epsilon_n+\epsilon_n^TH_n(s_{0n})\epsilon_n-\epsilon_n^TH_n(s) \epsilon_n.
\end{split}
\]
It is trivial that 
\[ 
 \epsilon_n^TH_n(s_{0n})\epsilon_n=p_{0n}(1+o_p(1)).  \ \ \mbox{(I)} 
\] 
We will show
\[
 \max\{\epsilon_n^TH_n(s)\epsilon_n,   \nu(s)  \leq
k_n \}= O_p( k_n \ln p_n),  \ \ \mbox{(II)}
\]
and 
\[ | \mu_n^T[ I_n-H_n(s)]\epsilon_n| = \sqrt{ \Delta_n(s)O_p(k_n \ln p_n ) },   \ \ \mbox{(III)}
\]
uniformly for all  $s$ with $\nu(s) \leq k_n$.   
Under the assumption of the theorem,  $2 k_n \ln p_n = o(n)$. Then,  by the asymptotic identifiability condition,  (I), (II) and (III) imply that 
\begin{equation}
\label{peq4}
y_n^T[I_n-H_n(s)]y_n-\epsilon_n^T[ I_n-H_n(s_{0n}) ]\epsilon_n \\
= \Delta_n(s)(1+o_p(1)),  
\end{equation}
uniformly for all  $s$ with $\nu(s) \leq k_n$.  
It then follows from (\ref{peq3}) and (\ref{peq4})  that 
\begin{equation}
\label{peq5}
 T_1 = n \ln \left( 1+ \frac{ \Delta_n(s) }{n} (1 + o_p(1)) \right), 
\end{equation}
uniformly for all $s$ with $\nu(s) \leq k_n$.

We now prove (II) and (III) in the following.
Let $ m = 2 k_n [\ln p_n   +  \ln (k_n \ln p_n)  ] $.  It  is obvious that $\dfrac{k_n}{m} \to 0$.  Note that we can express $\epsilon_n^TH_n(s)\epsilon_n = \chi^2_j(s) $ where $j=\nu(s)$.  By the Bonferroni inequality, we have 
\[ 
\begin{split}
& P(\max\{\epsilon_n^TH_n(s)\epsilon_n:  \nu(s) \leq k_n  \}\geq m ) \\
=  & P(\max\{\chi^2_j(s):  s \in \mathcal{S}_j, j  \leq k_n  \}\geq m )
\leq \sum\limits_{j=1}^{k_n}\tau(\mathcal{S}_j)P(\chi^2_j\geq m).
\end{split}
\]
By the fact that $\tau(\mathcal{S}_j) = {p_n \choose j} \leq p_n^j$ and Lemma 2, there is some $c$ close to 1, not depending on $j$ for $j \leq k_n$,  such that
\[
\begin{split}
             \tau(\mathcal{S}_j)P(\chi^2_j\geq m)  
\approx&  c\dfrac{1}{2^{j/2-1}\Gamma(j/2)}\dfrac{\tau(\mathcal{S}_j)  } {p_n^{k_n}}  
                 (k_n \ln p_n )^{- k_n}   m^{j/2-1} \\
          \leq & \frac{c}{m}    ( k_n \ln p_n )^{- j }   m^{j/2} 
          =   \frac{c}{m}   \left[ \sqrt {  \frac{m} {( k_n \ln p_n )^2 }  } \right]^{j} 
          =  \frac{c}{m}   q_n^j , \  \mbox{say}, 
\end{split}
\]
where 
 \[
     q_n  =   \sqrt {  \frac{m} {( k_n \ln p_n )^2 }  } 
             = \sqrt{  \frac{  2 [k_n \ln p_n  + k_n  \ln (k_n \ln p_n)  ] }  { ( k_n \ln p_n )^2 }  }(1+o(1))
             \leq q, 
\]
for some $q$  between 0 and 1, when $n$ is large enough, since $q_n \to 0$. Thus
\begin{equation}
\label{maxeq}
 P(\max\{\epsilon_n^TH_n(s)\epsilon_n:  \nu(s) \leq k_n  \}\geq m )
\leq  \frac{c}{m}   \sum\limits_{j=1}^{k_n}q^j  \leq  \frac{c}{m} \frac{q}{1-q} \to 0;
\end{equation} 
that is, 
\[ \max\{\epsilon_n^TH_n(s)\epsilon_n:  \nu(s) \leq k_n  \} =  m (1+o_p(1)) = O_p(k_n \ln p_n ) , \]
which establishes (II).

For verifying (III), note that we can express 
\[ \mu_n^T\{I_n-H_n(s)\}\epsilon_n=\sqrt{\Delta_n(s)}Z(s), \]
 where $Z(s) \sim N(0,1).$  For any $s$ with $\nu(s) \leq k_n$, we have 
\[
|\mu_n^T\{I_n-H_n(s)\}\epsilon_n| \leq \sqrt{\Delta_n(s)}\max \{|Z(s)|: \nu(s) \leq k_n \}.
\]
Let $m$ be the same as above.  Consider $P(\max \{|Z(s)|: \nu(s) \leq k_n \} \geq \sqrt{m})$. We have
\[
\begin{split}
     P( \max \{|Z(s)|: \nu(s) \leq k_n \} \geq \sqrt{m}) 
= & P(\max \{|Z(s)|:  s \in \mathcal{S}_j,  j  \leq k_n \} \geq \sqrt{m}) \\
\leq &  \sum_{j=1}^{k_n} \tau(\mathcal{S}_j) P(Z(s) \geq \sqrt{m})  
=  \sum_{j=1}^{k_n} \tau(\mathcal{S}_j) P(\chi_1^2 \geq m)   \\
\leq & \sum_{j=1}^{k_n} \tau(\mathcal{S}_j) P(\chi_j^2  \geq m), 
\end{split}
\]
since $P(\chi_1^2 \geq m)  < P(\chi_j^2  \geq m)$ by Lemma 2.  We have already shown that the last sum converges to zero.  This establishes (III). 

Now, putting (\ref{peq1}),  (\ref{peq2}) and (\ref{peq5}) together, we have
\[
\begin{split}
  & \mbox{{\rm EBIC}}_{\gamma}(s) - \mbox{{\rm EBIC}}_{\gamma}(s_{0n}) \\
=& n \ln \left( 1+ \frac{ \Delta_n(s) }{n} (1 + o_p(1)) \right) + (\nu(s) -p_{0n})\ln n+2\gamma(\ln \tau(\mathcal{S}_{\nu(s)} )-\ln \tau(\mathcal{S}_{p_{0n}})) \\ 
\geq & n \ln[1+ C p_{0n}\ln p_n/n(1 + o_p(1)) ] -p_{0n}(\ln n+2\gamma \ln p_n ) ,\\ 
\end{split}
\]
 for some positive $C$, when $n$ is large enough,  by the consistency condition.   Then by choosing $C > 1+2\gamma$, the difference goes to infinity as $n\to \infty$.

\vspace{0.15in}
\noindent {\it Case II:} $s_{0n} \subset s$.  
\vspace{0.15in}

When $ s_0 \subset s$, $\{I_n - H_n(s)\} X_n(s_0) = 0$.  Hence,
$ y_n^{T} \{ I_n - H_n(s)\} y_n = \epsilon_n^{T} \{ I_n - H_n(s) \} \epsilon_n $
and
\[
 \epsilon_n^{T} \{ I_n - H_n(s_0) \} \epsilon_n - \epsilon_n^{T} \{ I_n - H_n(s) \} \epsilon_n
=  \epsilon_n^{T} \{ H_n(s) - H_n(s_0)\} \epsilon_n
= \chi_{j}^2(s),
\]
where $\chi_{j}^2(s)$ is a $\chi^2$ random variable depending on $s$ with degrees of freedom $j$ and $j=\nu(s) - p_{0n}$. 
We obtain that
\begin{equation}
\label{peq6}
\begin{split}
       n \log \left(   \frac{ \epsilon_n^{T} \{ I - H_n(s_0) \} \epsilon_n  } {  \epsilon_n^{T} \{ I - H_n(s) \} \epsilon_n } \right) 
  =& n \log  \left \{ 
                    1 + \frac{\chi_{j}^2(s)}{\epsilon_n^{T} \{ I - H_n(s_0) \} \epsilon_n - \chi_{j}^2(s)}
                    \right \} \\
  \leq& \frac{n \chi_{j}^2(s)}{\epsilon_n^{T} \{ I - H_n(s_0) \} \epsilon_n - \chi_{j}^2(s)}. 
\end{split}
\end{equation}
As  $n \to \infty$, $n^{-1} \epsilon_n^{T} \{ I - H_n(s_0) \} \epsilon_n \to \sigma^2=1$, i.e., 
\begin{equation}
\label{peq7} 
\epsilon_n^{T} \{ I - H_n(s_0) \} \epsilon_n = n(1+o(1)).
\end{equation}  
Let $\tilde{\mathcal{S}}_j = \{ s: s \in \mathcal{S}_{j+p_{0n}}, s_0\subset s\}$.  Note that $\tau(\tilde{\mathcal{S}}_j ) = {p_n-p_{0n} \choose j} \leq p_n^j $.
Let $m_j = 2j [ \ln  p_n +  \ln ( j \ln p_n )]$.  In the same way as we derive (\ref{maxeq}), we have
\[
\begin{split}
       P(\max_{1 \leq j \leq k_n - p_{0n} } \dfrac{ \max\{ \chi_j^2 (s):  s\in  \tilde{\mathcal{S}}_j  \}}{m_j }\geq 1)
\leq &\sum_{j=1}^{k_n - p_{0n} } P(\max\{ \chi_j^2 (s):  s\in  \tilde{\mathcal{S}}_j  \}\geq m_j  )\\
\leq &\sum_{j=1}^{k_n - p_{0n} } \tau( \tilde{\mathcal{S}}_j ) P( \chi_j^2 \geq m_j  )
\leq \frac{1}{\ln p_n} \sum_{j=1}^{k_n - p_{0n} } q_j^j  \to  0,
\end{split}
\]
where 
\[ q_j = \sqrt{  \frac{ 2}{ j \ln p_n} + \frac{2 \ln (j \ln p_n) } { j (\ln p_n)^2}  } \leq \sqrt{  \frac{ 2}{  \ln p_n}(1 + o(1))   } \to 0.\]
Thus,
\begin{equation}
\label{peq8}
\max \{\chi_j^2(s): s  \in {\cal S}_{j+p_{0n}}, s_0 \subset s \}
= m_j \{1+ o_p(1)\}, 
\end{equation}
uniformly for all $s$ with $\nu(s) \leq k_n$ and $s_0 \subset s$.

It follows from (\ref{peq6}), (\ref{peq7}) and (\ref{peq8}) that 
\[
\begin{split}
        n \log \left(   \frac{ \epsilon_n^{T} \{ I - H_n(s_0) \} \epsilon_n  } {  \epsilon_n^{T} \{ I - H_n(s) \} \epsilon_n } \right)  
    \leq  & \frac{n m_j }{ [n - m_j (1+o_p(1))] } \\
     \leq  & m_j (1+o_p(1)) \leq 2 j (1+\delta) \ln  p_n (1+o_p(1)), 
\end{split}
\]
uniformly for all $s$ with  $\nu(s) \leq k_n$ and $s_0 \subset s$,  noting that $m_j \leq 2j [ \ln p_n + \ln ((k_n - p_{0n}) \ln p_n)]  = 2j  (1+\delta)  \ln  p_n (1+o_p(1)) $ and $ m_{j} = 2j (1+\delta)  \ln  p_n (1+o_p(1)) $ for $j=k_n - p_{0n}$.  Thus
\[ T_1 \geq - 2 j (1+\delta) \ln  p_n (1+o_p(1)).
\]
When $p_{0n} \leq \nu(s) \leq k_{n}$ we have $\ln \nu(s) / \ln p_n \to \delta $ uniformly,  hence, by Lemma 1, 
\[ T_2 = j \ln n+2\gamma (1-\delta) j \ln p_n(1+o(1)).
\]
 Finally we have 
\[
\begin{split}
  & \mbox{{\rm EBIC}}_{\gamma}(s) - \mbox{{\rm EBIC}}_{\gamma}(s_{0n}) \\
\geq & j \ln n+2\gamma (1-\delta) j \ln p_n(1+o(1))- 2 j (1+\delta) \ln  p_n   (1+o_p(1)) > 0, 
\end{split}
\]
uniformly for all $s$ with  $\nu(s) \leq k_n$ and $s_0 \subset s$, if $n$ is big enough, when $\gamma > \frac{1+\delta}{1-\delta} - \frac{\ln n }{ 2 (1-\delta) \ln p_n }$.

\end{proof}

\begin{center}
\begin{table}
\caption{ The PDR and FDR of the SIS-SCAD-EBIC procedure with
Structure I (power decay correlation) averaged over 200 replicates
(the numbers in parentheses are standard deviations) }
\label{table1}

\begin{tabular}{cc|ccc|ccc} \hline \hline
\multicolumn{8}{c}{$c=1$} \\ \hline
 & & \multicolumn{3}{|c}{PDR} &  \multicolumn{3}{|c}{FDR} \\ \cline{3-8}
$n$ & $h$ &   $\gamma_{\mbox{\tiny BIC}}$ & $\gamma_{\mbox{\tiny SC}} $ &
$\gamma_{\mbox{{\tiny  mBIC}}}$  &  $\gamma_{\mbox{\tiny BIC}}$ & $\gamma_{\mbox{\tiny SC}} $ &
$\gamma_{\mbox{{\tiny  mBIC}}}$ \\ \hline
100 & .4 &.726(.242) & .450(.291) &.384(.288) &.571(.212) &.074(.205) &.050(.181) \\
       & .6 &.861(.187)  & .700(.271)  &.633(.301)  &.478(.216)  & .080(.170)  &.044(.123)  \\
       & .8 &.973(.090)   & .921(.159) & .909(.176)  &.363(.204)  &.085(147) & .056(.120) \\ \hline
 200 & .4  &.759(.205)   & .532(.270)  &.467(.270)  &.662(.177)   & .034(.101) & .017(.071) \\
       & .6 &.910(.144)  & .758(.256) &.711(.282)  &.574(.185)   & .080(.145)  &.038(.100)  \\
       & .8 &.989(.056)  & .957(.105) & .947(.128)  &.389(.200) & .060(.115)  &.045(.105) \\  \hline
500 & .4 &.826(.146)   & .640(.212)  &.604(.214)  &.768(.100)   & .037(.090)  &.011(.046) \\
       & .6  &.943(.100)  & .863(.164)  &.836(.181)  &.660(.133)   & .066(.128)  &.028(.079) \\
       & .8 &.994(.035)   & .983(.060)  &.980(.067)  &.504(.190)  & .027(.073)  &.019(.065) \\ \hline
1000 & .4 &1.000(.00)  & .999(.008)  & .999(.011)  &.662(.024)  & .019(.041)  &.009(.028)  \\
       & .6 &1.000(.00)   & 1.000(.00)  & 1.000(.00)  &.531(.037)  & .019(.041)  & .008(.026)\\
       & .8 &  1.000(.00) &  1.000(.00)& 1.000(.00) &.470(.010)   & .007(.025)  & .002(.014)\\  \hline
\multicolumn{8}{c}{$c=2$} \\ \hline
 & & \multicolumn{3}{|c}{PDR} &  \multicolumn{3}{|c}{FDR} \\ \cline{3-8}
$n$ & $h$ &   $\gamma_{\mbox{\tiny BIC}}$ & $\gamma_{\mbox{\tiny SC}} $ &
$\gamma_{\mbox{{\tiny  mBIC}}}$  &  $\gamma_{\mbox{\tiny BIC}}$ & $\gamma_{\mbox{\tiny SC}} $ &
$\gamma_{\mbox{{\tiny  mBIC}}}$ \\ \hline
100 & .4 &  .531(.183)&.243(.169) &.198(.162)   &  .507(.222) &.069(.204) & .041(.172)  \\
       & .6 &  .680(.166) &.416(.213) &.350(.206)  &  .447(.187) & .074(.173) &.026(.093) \\
       & .8 &  .850(.153) &.708(.225) &.628(.248)   & .373(.163) &.118(.143) &.068(.118) \\ \hline
200 & .4 & .613(.162) &.306(.164)  &.260(.161)  & .619(.162) &.028(.096) &.010(.066)  \\
       &.6 &.720(.148)&.518(.211)  &.456(.207) & .545(.181)&.036(.082) & .018(.061)  \\
       & .8 & .895(.125) &.745(.199) &.703(.217) & .447(.164) &.086(.117) &.053(.096)  \\ \hline
500 & .4 &  .732(.130) &.425(.174) &.371(.166)  &  .774(.076) &.014(.054) & .004(.025) \\
       & .6 & .832(.104)  &.635(.176)  &.590(.186)  &  .695(.112)& .028(.064) & .010(.031)  \\
       & .8 &.956(.067) &.875(.135)  &.847(.157)  & .535(.159)&.098(.121) & .068(.104)  \\ \hline
1000   & .4 & .758(.108)  &.537(.161)  &.491(.164) &   .825(.055) & .012(.040) & .005(.025) \\
         &   .6   &  .849(.102)  &.715(.134)  &.689(.144)  & .761(.077) &.025(.062) &.010(.035)  \\
        &  .8    & .969(.054) &.925(.084) &.906(.106)  &.581(.146) &.095(.110)&.072(.095) \\ \hline
\end{tabular}
\end{table}
\end{center}


\begin{center}
\begin{table}[h]
\caption{ The PDR and FDR of the SIS-SCAD-EBIC procedure with
Structure II (block with equal pairwise correlation) averaged over 200
replicates (the numbers in parentheses are standard deviations) }
\label{table2}

\begin{tabular}{cc|ccc|ccc} \hline \hline
\multicolumn{8}{c}{$c=1$} \\ \hline
 & & \multicolumn{3}{|c}{PDR} &  \multicolumn{3}{|c}{FDR} \\ \cline{3-8}
$n$ & $h$ &   $\gamma_{\mbox{\tiny BIC}}$ & $\gamma_{\mbox{\tiny SC}} $ &
$\gamma_{\mbox{{\tiny  mBIC}}}$  &  $\gamma_{\mbox{\tiny BIC}}$ & $\gamma_{\mbox{\tiny SC}} $ &
$\gamma_{\mbox{{\tiny  mBIC}}}$ \\ \hline
 100 & .4     & .733(.285)  & .402(.318) & .343(.291)  & .427(.268)   & .229(.369) &.198(.362)  \\
       & .6  & .933(.154)  & .772(.297)  &.703(.321) & .340(.213)  & .117(.197) & .094(.207)  \\
       & .8 &.996(.042)   & .967(.118) & .960(.125)   & .293(.203) & .053(.132)  &.036(.114)  \\ \hline
 200 & .4 & .868(.203)  & .534(.303)  &.479(.306) &.442(.206)  &.133(249)  &.110(.246)  \\
       & .6 & .994(.039)   & .931(.168)  &.889(.214)   & .321(.173) & .107(161)  &.078(.143)  \\
       & .8 & 1.000(.00) & .996(.031) & .994(.040)  &.292(.165)  & .025(.081) & .017(.070)  \\ \hline
500 & .4  & .948(.093)   & .754(.178)  &.723(.184)  & .689(.114)  & .056(.107) &.049(.103) \\
       & .6 & .993(.035)  & .922(.121)  &.904(.132)    &.626(.127) & .031(.080) & .019(.064)  \\
       & .8 & 1.000(.00)  & .997(.024)  & .992(.044)    & .585(.151)  & .060(.110)  & .031(.083)  \\ \hline
1000 & .4 &.940(.080)   & .813(.158)  &.785(.180)& .818(.046)  &  .073(.113) &  .049(.092) \\
       & .6 &.995(.025)& .988(.041) &.986(.043) & .739(.066)  & .039(.084)  &.035(.079) \\ 
       & .8 &.999(.010)  & .998(.017)  &.996(.024) &.653(.107)   & .024(.070)  &.017(.061) \\ \hline
\multicolumn{8}{c}{$c=2$} \\ \hline
 & & \multicolumn{3}{|c}{PDR} &  \multicolumn{3}{|c}{FDR} \\ \cline{3-8}
$n$ & $h$ &   $\gamma_{\mbox{\tiny BIC}}$ & $\gamma_{\mbox{\tiny SC}} $ &
$\gamma_{\mbox{{\tiny  mBIC}}}$  &  $\gamma_{\mbox{\tiny BIC}}$ & $\gamma_{\mbox{\tiny SC}} $ &
$\gamma_{\mbox{{\tiny  mBIC}}}$ \\ \hline
100 & .4 &  .430(.239)&.193(.174) &.173(.164)  &.449(.294) &.310(.411) &.295(.408) \\
       & .6 &  .684(.234) &.390(.236)  &.343(.224) &.343(.220)& .164(.235)&.150(.253)  \\
       & .8 & .881(.179) &.676(.266)  &.603(.284)   &.308(.194) & .105(.174) &.096(.175)  \\ \hline
200 & .4 &  .489(.206)  &.199(.142) &.165(.133)  &.416(.235) &.134(.275) &.115(.259)  \\
       & .6 &  .727(.192) &.421(.227) &.356(.214)   & .351(.195) &.065(.144)&.055(.132) \\
       & .8 &  .919(.135)&.718(.254) &.672(.269)   & .351(.184) &.055(.099) &.043(.088) \\  \hline
500 & .4 &  .664(.137)&.258(.132) &.238(.132)  &.669(.145) & .031(.099) & .020(.076) \\
       & .6 & .834).127) &.468(.211) &.407(.209) &  .609(.132) &.029(.073) & .014(.047)  \\
       & .8 &  .944(.094) &.804(.244) &.778(.266)   &  .485(.198)& .084(.108) & .068(.095) \\ \hline
1000 & .4 & .675(.133)&.311(.158)  &.284(.158) &.830(.079) &.017(.055)&.014(.050)  \\
       & .6 &.882(.134) &.551(.234)  &.496(.240)  & .744(.115) &.060(.108) & .033(.073)  \\
       & .8 &.960(.078)&.884(.195) &.877(.202)&.616(.178) &.069(.099)&.061(.087)  \\ \hline
\end{tabular}
\end{table}
\end{center}


\begin{center}
\begin{table}
\caption{ The PDR and FDR of the SIS-SCAD-EBIC procedure with
Structure III (block with uniformly distributed eigenvalues) averaged
over 200 replicates (the numbers in parentheses are standard
deviations) } 
\label{table3}

\begin{tabular}{cc|ccc|ccc} \hline \hline
\multicolumn{8}{c}{$c=1$} \\ \hline
 & & \multicolumn{3}{|c}{PDR} &  \multicolumn{3}{|c}{FDR} \\ \cline{3-8}
$n$ & $h$ &   $\gamma_{\mbox{\tiny BIC}}$ & $\gamma_{\mbox{\tiny SC}} $ &
$\gamma_{\mbox{{\tiny  mBIC}}}$  &  $\gamma_{\mbox{\tiny BIC}}$ & $\gamma_{\mbox{\tiny SC}} $ &
$\gamma_{\mbox{{\tiny  mBIC}}}$ \\ \hline
100 & .4 &  .915(.146)& .667(.302)  &.564(.327)   & .428(.191)  & .041(.102)& .020(.078) \\
       & .6 &  .996(.031)  & .964(.116)  &.950(.133)  & .360(.181)   & .046(.105)  &.019(.063) \\
       & .8 & 1.000(.00)& 1.000(.00)  & 1.000(.00)   & .326(.165)  & .038(.096)  & .011(.051)  \\ \hline
 200 & .4 & .993(.037) & .865(.206) &.811(.252) & .575(.162)  &.050(.101) &.024(.073)  \\
       & .6 & 1.000(.00)  & .999(.014)  & .999(.014)  & .536(.129)  & .032(.081) & .013(.048)  \\
       & .8 & 1.000(.00)   & 1.000(.00) & 1.000(.00)  & .457(.138)  & .023(.065)  & .009(.042)  \\ \hline
500 & .4 &  1.000(.00)  & .971(.081)  & .961(.090)  &  .768(.042)  & .041(.075)  &.023(.055)  \\
       & .6 & 1.000(.00)  &1.000(.00)  & 1.000(.00) & .704(.058)   & .022(.060) & .010(.043) \\
       & .8 & 1.000(.00)   & 1.000(.00)  & 1.000(.00)  & .608(.091)   & .016(.050)  & .007(.038)  \\ \hline
1000 & .4 &1.000(.00) & .999(.011)  & .997(.017)   &.790(.040) & .023(.046)  &.008(.028)\\
       & .6 &1.000(.00)   & 1.000(.00)& 1.000(.00)  &.740(.038)  & .018(.041) & .005(.021) \\
       & .8 &1.000(.00)   & 1.000(.00)  & 1.000(.00) &.705(.051)   & .005(.022) & .002(.012) \\ \hline
\multicolumn{8}{c}{$c=2$} \\ \hline
 & & \multicolumn{3}{|c}{PDR} &  \multicolumn{3}{|c}{FDR} \\ \cline{3-8}
$n$ & $h$ &   $\gamma_{\mbox{\tiny BIC}}$ & $\gamma_{\mbox{\tiny SC}} $ &
$\gamma_{\mbox{{\tiny  mBIC}}}$  &  $\gamma_{\mbox{\tiny BIC}}$ & $\gamma_{\mbox{\tiny SC}} $ &
$\gamma_{\mbox{{\tiny  mBIC}}}$ \\ \hline
100 & .4 &  .643(.218)  &.240(.201) &.155(.179)  & .409(.206) &.071(.185) &.028(.128) \\
       & .6 &  .911(.141) &.589(.298)  &.461(.302) & .346(.168) & .092(.163)&.045(129)  \\
       & .8 &  .995(.033)  &.975(.100)&.964(.135) &.237(.136)&.089(.101)&.069(.092)  \\ \hline
200 & .4 &  .801(.147) &.307(.210)  &.209(.179) &  .536(.136)&.050(.142) &.013(.061) \\
       & .6 &  .974(.063)&.817(.198)  &.742(.236)  &  .443(.147) &.076(.095)&.045(.073)  \\
       & .8 &.999(.010) &.993(.041)  &.989(.048) & .322(.121) & .046(.074) & .034(.063) \\ \hline
500 & .4 &  .933(.076) &.578(.204)&.451(.215)   &  .723(.079)& .035(.067)& .009(.036) \\
       & .6 &.992(.030)  &.946(.073) &.930(.094)  &  .642(.091) & .062(.078) & .045(.069) \\
       & .8 &  .999(.005)  &.998(.016) &.997(.017)  &  .498(.105) & .023(.044) &.014(.036) \\ \hline
1000 & .4 &.970(.049) &.780(.170) &.688(.207)& .809(.051) &.042(.063) & .018(.039) \\
       & .6 &.997(.013) &.976(.054)  &.973(.058)  & .738(.059) &.030(.053)& .024(.042)  \\
       & .8 & .999(.004) &.998(.011)  &.998(.012) &.608(.085) &.011(.031) & .006(.022)  \\ \hline
\end{tabular}
\end{table}
\end{center}

\end{document}